\documentclass[10pt,a4paper]{article}
\usepackage{amssymb, amsmath,amsfonts,amsthm}
\usepackage{graphicx}
\usepackage{color}
\usepackage[all]{xy}
\usepackage{psfrag}

\theoremstyle{plain}

\theoremstyle{definition}

\newcommand{\be}{\begin{enumerate}}
\newcommand{\ee}{\end{enumerate}}

\begin{document}
\bibliographystyle{plain}

\author{S. O Rourke}
\title{A tree-free group that is not orderable}\maketitle

\begin{abstract}
I.~M. Chiswell has asked whether every group that admits a free
isometric action (without inversions) on a $\Lambda$-tree is
orderable. We give an example of a multiple HNN extension $\Gamma$
which acts freely on a $\mathbb{Z}^2$-tree but which has non-trivial
generalised torsion elements. The existence of such elements implies
that $\Gamma$ is not orderable.\footnote{I would like to thank Ian Chiswell for helpful conversations.}
\end{abstract}

\maketitle

Let $\Lambda$ be an ordered abelian group. A group is
\emph{$\Lambda$-free} if it admits a free isometric action without
inversions on a $\Lambda$-tree, and \emph{tree-free} if it is
$\Lambda'$-free for some $\Lambda'$. We refer to the book
\cite{Chiswell-book} for a detailed account of the fundamentals of
$\Lambda$-trees.

In this book Chiswell asks \cite[\S5.5 Question 3]{Chiswell-book}
whether all tree-free groups are orderable, or at least
right-orderable. There has recently been some progress made on
questions of orderability in tree-free groups. Chiswell himself has
shown \cite[Theorem 3.8]{Chiswell-ro} that $\mathbb{R}^n$-free
groups are right-orderable. Kharlampovich, Myasnikov and Serbin have
shown \cite[Corollary 4]{KM-Lambda} that finitely presented
tree-free groups are $\mathbb{R}^n$-free for some $n$; thus these
groups are right-orderable. Chiswell has shown moreover
\cite[Theorem 4.5]{Chiswell-LIO} that tree-free groups admit a
locally invariant order.

In their recent survey Kharlampovich, Myasnikov and Serbin state
\cite[Corollary 19]{KMS-long-survey} that finitely presented
tree-free groups have a finite index subgroup that embeds in a
right-angled Artin group: this
is a consequence of their result \cite[Theorem 2]{KM-Lambda} and the
extensive work of Wise (see \cite[\S16]{Wise-hierarchy} and
\cite{Wise-shorter-hierarchy}) on quasi-convex hierarchies on
groups.
Since right-angled Artin groups are residually torsion-free
nilpotent (see \cite[Chapter 3, Theorem 1.1]{Droms-thesis}), it follows that finitely presented
tree-free groups are virtually residually torsion-free nilpotent;
hence they are \emph{virtually} orderable.

The author has recently \cite{affine-paper} raised the question of
whether $\mathbb{Z}^n$-free groups are residually torsion-free
nilpotent. An affirmative answer to this question would have implied
that $\mathbb{Z}^n$-free groups are orderable.\medskip

Nevertheless the answer to Chiswell's question is negative, even
when restricted to finitely presented $\mathbb{Z}^2$-free groups, as
will show presently. It follows that the word `virtually' cannot be
dropped in the discussion above. This suggests an analogy with the
situation of braid groups $B_n$ and their finite index subgroups,
the pure braid groups $P_n$: the former are right-orderable but not
orderable (see \cite[\S4]{Rolfsen-braid-survey}), while the latter
are residually torsion-free nilpotent \cite{Falk-Randell}.

Recall that a group $G$ is \emph{orderable} if there is a linear
order $\leq$ on $G$ satisfying $x\leq y\Rightarrow gxh\leq gyh$ for
$g,h\in G$. (One can define right-orderable by restricting to $g=1$
in the definition above.) It is well-known and easy to see that in
an orderable group $G$ there can be no non-trivial \emph{generalised
torsion elements}: these are elements $g$ such that
$g^{h_1}g^{h_2}\cdots g^{h_n}=1$ for some $h_1,\ldots,h_n\in G$ and
$n\geq 1$. (Here $g^h$ denotes the conjugate $h^{-1}gh$.)

Let $F$ be the free group on $\{x,y,z\}$, and consider the natural
free action of $F$ on the corresponding Cayley graph, viewed as a
$\mathbb{Z}$-tree. Observe that $xy^{-1}$, $yz^{-1}$ and $zx^{-1}$
and their respective inverses belong to distinct conjugacy classes
since they are cyclically reduced as elements of $F$ and none is a
cyclic permutation of another. Moreover, the translation lengths of
these elements are all equal to 2.

Now taking $s_1=xy^{-1}=s_2$, $t_1=yz^{-1}$, $t_2=zx^{-1}$, $u=u_1$
and $v=u_2$, and applying \cite[Proposition 4.19]{Bass}, the
multiple HNN extension \begin{center}$\begin{array}{rl}\Gamma
&=\langle x,y,z,u_1,u_2\ |\ u_1s_1u_1^{-1}=t_1,\
u_2s_2u_2^{-1}=t_2\rangle\\ &=\langle u,v,x,y,z\ |\
u(xy^{-1})u^{-1}=yz^{-1},\
v(xy^{-1})v^{-1}=zx^{-1}\rangle\end{array} $\end{center} is seen to
be $\mathbb{Z}^2$-free. However,
$$1=(xy^{-1})(yz^{-1})(zx^{-1})=xy^{-1}\cdot u(xy^{-1})u^{-1}\cdot
v(xy^{-1})v^{-1},$$ whence $xy^{-1}$ is a non-trivial generalised
torsion element of $\Gamma$, and $\Gamma$ is not orderable. This
gives the promised negative answer to Chiswell's question.

\def\cprime{$'$}
\providecommand{\bysame}{\leavevmode\hbox
to3em{\hrulefill}\thinspace}
\providecommand{\MR}{\relax\ifhmode\unskip\space\fi MR }
\providecommand{\MRhref}[2]{%
  \href{http://www.ams.org/mathscinet-getitem?mr=#1}{#2}
} \providecommand{\href}[2]{#2}

\vspace{0.5cm}

\begin{minipage}[t]{3 in}
\noindent Shane O Rourke\\ Department of Mathematics\\
Cork Institute of Technology\\
Rossa Avenue\\ Cork\\ IRELAND
\\ \verb"shane.orourke@cit.ie"
\end{minipage}


\begin{thebibliography}{10}

\bibitem{Bass}
H. Bass, \emph{Group actions on non-{A}rchimedean trees}, Arboreal
group
  theory ({B}erkeley, {CA}, 1988), Math. Sci. Res. Inst. Publ., vol.~19,
  Springer, New York, 1991, pp.~69--131. \MR{1105330 (93d:57003)}

\bibitem{Chiswell-book}
I.~M. Chiswell, \emph{Introduction to {$\Lambda$}-trees}, World
Scientific
  Publishing Co. Inc., River Edge, NJ, 2001. \MR{1851337 (2003e:20029)}


\bibitem{Chiswell-LIO}
\bysame, \emph{Locally invariant orders on groups}, Internat. J.
Algebra
  Comput. \textbf{16} (2006), no.~6, 1161--1179. \MR{2286427 (2007k:20083)}

\bibitem{Chiswell-ro}
\bysame, \emph{Right orderability and graphs of groups}, J. Group
Theory
  \textbf{14} (2011), no.~4, 589--601. \MR{2818951}


\bibitem{Droms-thesis}
C.~Droms, \emph{Graph groups}, Ph.D. thesis, Syracuse University,
1983, Available at http://educ.jmu.edu/~dromscg/vita/thesis/.

\bibitem{Falk-Randell}
M. Falk and R. Randell, \emph{Pure braid groups and products of free
  groups}, Braids ({S}anta {C}ruz, {CA}, 1986), Contemp. Math., vol.~78, Amer.
  Math. Soc., Providence, RI, 1988, pp.~217--228. \MR{975081 (90d:20070)}

\bibitem{KM-Lambda}
O.~Kharlampovich, A.~Myasnikov, and D.~Serbin, \emph{Groups acting
freely on
  ${\Lambda}$-trees}, (arXiv:0911.0209v4), November 2011.

\bibitem{KMS-long-survey}
\bysame, \emph{Actions, length functions, and non-archimedean
words},
  (arXiv:1211.3207v1), November 2012.

\bibitem{affine-paper}
S.~O~Rourke, \emph{Affine actions on non-archimedean trees},
Internat.
  J. Algebra Comput. (to appear), (arXiv:1112.4832v2).

\bibitem{Rolfsen-braid-survey}
D. Rolfsen, \emph{New developments in the theory of {A}rtin's braid
groups},
  Proceedings of the {P}acific {I}nstitute for the {M}athematical {S}ciences
  {W}orkshop ``{I}nvariants of {T}hree-{M}anifolds'' ({C}algary, {AB}, 1999),
  vol. 127, 2003, pp.~77--90. \MR{1953321 (2004f:20070)}

\bibitem{Wise-hierarchy}
D.~Wise, \emph{The structure of groups with a quasiconvex
hierarchy}, Available
  at http://www.math.mcgill.ca/wise/papers.html, 2011.

\bibitem{Wise-shorter-hierarchy}
\bysame, \emph{Research announcement: the structure of groups with a
  quasiconvex hierarchy}, Electron. Res. Announc. Math. Sci. \textbf{16}
  (2009), 44--55. \MR{2558631 (2011c:20052)}

\end{thebibliography}
\end{document}